\documentclass[12pt]{article}
\usepackage{amsmath}
\usepackage{amsfonts}
\usepackage{amssymb}

\newtheorem{thm}{Teorema}
\numberwithin{thm}{section}
\newtheorem{Lemma} [thm]{{\textsc Lemma}}

\newtheorem{Prop} [thm] {\textsc{Proposition}}

\newcommand{\QED}{\hfill $\square$\bigskip}
\newcommand{\pf}{{\noindent\bf Proof.\ \ }}

\newcommand{\semidir}{\kern-2pt>\kern-6pt\triangleleft}
\newcommand{\Aut}{\mbox{\rm Aut}}
\newcommand{\PAut}{\mbox{\rm PAut}}
\newcommand{\IAut}{\mbox{\rm IAut}}
\newcommand{\Fit}{\mbox{\rm Fit}}
\newcommand{\Z}{{\Bbb Z}}
\newcommand{\Q}{{\Bbb Q}}

\newcommand{\la}{{\langle}}
\newcommand{\ra}{{\rangle}}
\newcommand{\g}{{\gamma}}%
\newcommand{\G}{{\Gamma}}%

\newcommand{\TI}{{${\it \tilde T}$}}
\newcommand{\iso}{\simeq}

\newcommand{\sn}{{\rm \kern1pt sn \kern1pt}}
\newcommand{\an}{{\rm \kern2pt an \kern1pt}}
\newcommand{\cf}{{\rm \kern2pt cf \kern1pt}}
\newcommand{\nf}{{\rm \kern2pt nf \kern1pt}}
\newcommand{\f}{{\rm \kern2pt f \kern1pt }}
\newcommand{\ff}{{\kern1pt \rm \bar{f}\kern1pt }}
\newcommand{\C}{\rm \kern2pt c \kern1pt}
\newcommand{\cn}{\rm \kern2pt cn\kern1pt}
\newcommand{\ssim}{\rm \kern2pt \approx}

\begin{document}


\title{\vskip-2cm On soluble groups whose\\ subnormal subgroups are inert
}

\author{Ulderico Dardano\footnote{Dipartimento di Matematica e
Applicazioni ``R.Caccioppoli'', Universit\`a di Napoli ``Federico
II'', Via Cintia - Monte S. Angelo, I-80126 Napoli, Italy.
\texttt{dardano@unina.it}}\  \ -\
Silvana Rinauro\footnote{
Dipartimento di Matematica, Informatica ed Economia, Universit\`a della
Basilicata, Viale dell'Ateneo Lucano 10,
I-85100 Potenza, Italy. \texttt{silvana.rinauro@unibas.it} } }
\date{}

\maketitle

\noindent{\textbf{Abstract.}  A subgroup H of a group G is called inert if for each $g\in G$ the index of $H\cap H^g$ in $H$ is finite. We  give a classification
of soluble-by-finite groups $G$ in which subnormal subgroups are inert in the cases where $G$  has no nontrivial torsion normal subgroups or $G$ 
 is finitely generated.
}
\bigskip




\section{Introduction and Main Results}

According to \cite{BKS} and \cite{R},  a subgroup H of a group G is said to be {\it inert} if for each $g\in G$ the index of $H\cap H^g$ in $H$ is finite. This is equivalent to saying that $H$ is commensurable with each of its conjugates. Recall that two subgroups $H$ and $K$ of a group are said \emph{commensurable} if and only if both indices $|H:(H\cap K)|$
and $|K:(H\cap K)|$ are finite. 

A group all of whose subgroups are inert is called {\it inertial} in \cite{R} where, in the context of generalized soluble groups  (with some finiteness conditions), a characterization of inertial groups was given.  Inertial groups received attention also  in the context of locally finite groups (see \cite{BKS} and \cite {DET} for example).

Recall that  Theorem A  of \cite{R}
 states that  \emph{a hyper-(abelian or finite)  group $G$ without non-trivial torsion normal subgroups is inertial if and only if  $G$ is abelian or dihedral.}

In Theorem B of \cite{R} it is shown that \emph{a finitely generated hyper-(abelian or finite) group $G$ is inertial
if and only if it has a torsion-free abelian normal subgroup of finite index in which elements of $G$ induce
power automorphisms}. Actually, this is equivalent to the fact that $G$ has a finite normal   subgroup $F$ such that $G/F$ is torsion-free abelian or dihedral on a  torsion-free abelian group (see Proposition 6.1 of \cite{R}).

\medskip

For terminology, notation and basic facts we refer to \cite{LR} and \cite{Rm}.
In particular, a {\it dihedral group} $G$ on an abelian group $A$ is a 
  group $ G=\la x\ra\ltimes A$ where $x$ acts faithfully on $A$ as the inversion map. 
Also, an  automorphism $\g$ of a group $A$ is said a {\it power automorphism} if and only if $\g$ maps each subgroup into itself. Thus $a^\g=a^n$ for all $a\in A$ where $n = n(a)\in \Z$. Clearly, if $A$ is non-torsion, then $n=\pm1$ . 

\medskip

The class {\it T}\ of groups in which subnormal subgroups are normal and its generalizations
received much attention in the literature.
In the present work we consider the class \TI \ of groups whose {\it subnormal subgroups are inert}. 

In Proposition \ref{PROPin+sn->sin} we show that  \emph{an inert subnormal subgroup $H$ of a group $G$ is strongly inert} (according to the terminology of  \cite{DFMT}), 
that is $H$ has the property that $|\la H,H^g\ra:H|$ is finite for all $g\in G$. Clearly, strongly inert subgroups are inert. Note that {in a dihedral group $G$ each subgroup $H$ is  inert}, since $|H:H_{G}|\le 2$. However, if $G$ is infinite, then non-normal 
subgroups of order $2$ are not strongly inert.
Recall that  $H^G$  (resp. $H_G$) denotes the smallest (resp. largest) normal subgroup of $G$ containing $H$ (resp. contained in $H$).  
Groups all of whose subgroups are strongly inert have been studied in \cite{DFMT}. 

Before stating our main results, we introduce some terminology. An automorphism of a (possibly non-abelian) group $A$ is said to be an {\em inertial automorphism} if and only if it maps each subgroup of $A$ to a commensurable one. Clearly a \TI-group acts on its normal abelian sections by means of inertial automorphisms. 
In \cite{DR1} we noticed that  
inertial automorphisms of a group $A$ form a group $\IAut(A)$ and that, {\it if $A$ is torsion-free abelian, then  $\IAut(A)$ consists of maps $\g$ such that there are coprime integers $m,n$ (with $n>0$) } such that
$$(a^\g)^{n}=a^{m}\ \forall a\in A$$ 
where clearly $m,n$ are uniquely defined and {\it $mn=\pm1$ if $A$ has infinite rank} (see also Proposition \ref{Recall}(1) below).  Recall a group is said to have finite (Pr\"ufer) rank if there is an integer $r$ such that 
every finitely generated subgroup can be generated by $r$ elements.

\newpage

 For a complete description of inertial automorphisms of any abelian group $A$ see \cite{DR2} and \cite{DR4}.

 \medskip
 \noindent{\bf Definition} \emph{
  A group $G$ is said to be  {\em semidihedral} on  a torsion-free abelian subgroup $A\ne 1$ if $G$ acts on $A$ by means of inertial automorphims and $C_{G}(A)= A$.}
\bigskip

 The reader is warned that the word {\em semidihedral} has been  also used  with a different meaning in other areas of group theory.  Note that, by Proposition \ref{Recall}(1) below,  it follows that in the  above conditions we have that $G/A$ is abelian. Moreover, $A=\Fit(G)$ is uniquely determined, since the elements of $G\setminus A$ act on $A$ as fixed-point-free automorphisms. Furthermore, if elements of $G$ induce on $A$ power automorphisms (that is the identity and the inversion map, since $A$ is torsion-free), then $G$ is abelian or dihedral.  This is the case when $A$ has infinite rank. 

We are now in a position to state a theorem corresponding  to Theorem A of \cite{R}.  Notice  that clearly \emph{a semidihedral group has no non-trivial torsion normal subgroups}.

\medskip \noindent{\bf Theorem $\bf \tilde A$} \emph{A hyper-(abelian or finite) group $G$ without non-trivial torsion normal subgroups  is a \TI-group if and only if $G$ is semidihedral on a torsion-free abelian subgroup.}

\medskip


Then we have a statement that corresponds to Theorem B in \cite{R} and answers  Question C in \cite{dG}.

\bigskip
\noindent{\bf Theorem $\bf \tilde B$} \emph{For a  finitely generated soluble-by-finite group $G$, the following are equivalent:\\
$i)$\ \ $G$ is a \TI-group;\\
\smallskip
$ii)$\ $G$ has a finite normal subgroup $F$ such that $G/F$ is a semidihedral group.} 

Actually, for a finitely generated group $G$, condition $(ii)$ is equivalent to the fact that $G$ has a semidihedral normal subgroup with finite index $G_0$ such that $G$ acts by means of inerial automorphisms on $Fit(G_{0})$ and trivially on $G/G_0$ (see Proposition \ref{Prop2SuSemidiedrali} for details)

\newpage

\section{\hskip-5mm Preliminary results and proof of Theorem $\bf \tilde A$}

\medskip

Our first result seems to be missing in the literature. 

\begin{Prop}\label{PROPin+sn->sin}
Inert subnormal subgroups are strongly inert.
\end{Prop}

\medskip
\pf This follows from next lemma.

\medskip

\begin{Lemma}\label{LEMMAin+sn->sin}
Let $H$ be an inert subnormal subgroup of a group $G$ and $K\le G$. If\ $|K:(H\cap K)|$ is finite,  then $|\la H,K\ra:H|$ is finite as well.
\end{Lemma}

\pf We may assume $G=\la H,K\ra=H^GK$ and 
proceed by induction on the subnormality defect $d$ of $H$ in $G$, the statement being trivial if $d\le 1$.
 Since $|G:H|=|G:H^G|\cdot|H^G:H|$
and $|G:H^G|\le |K:(H\cap K)|=:n$ is finite, we only have to show that 
$|H^G:H|$ is finite. To this aim note that $H^G$ is the join of at most $n$ conjugates of $H$. Thus the statement follows from the following claim: \emph{for any positive integer $r$, any subgroup $H^+$  generated by at most $s$ conjugates of $H$ is commensurable with $H$}. To show this fact proceed  by induction on $s$, the claim being trivial if $s=1$ since $H$ is inert by hypothesis. Assume that the claim is  true for $s$ and consider a subgroup $H^+=\la H_1,\ldots,H_{s+1}\ra$ where each $H_i$ is conjugate to $H$. Denote commensurability by $\sim$ and note that it is a transitive relation. Then, by induction on $s$, we have that $L:=\la H_1,\ldots,H_{s}\ra\sim H\sim H_{s+1}$.
Further, since  $H_{s+1}$ has subnormality defect at most $d-1$ in $H^G$, we may apply induction (on $d$)  to the group $H^{G}$ and its subgroups $H_{s+1}$ and $L$. We have that $|\la H_1,\ldots,H_{s+1}\ra:H_{s+1}|$ is finite. Thus $H^+\sim H_{s+1}\sim H$. \QED

Recall that if an abelian group $A$ is bounded, the group $\PAut(A)$ of power automorphism of $A$ is well-known to be finite. If $A$ is non-torsion abelian, then the only power automorphisms of $A$ are the identity and the inversion map (see \cite{Rm}).

\smallskip

By next proposition we recall some facts about the group  $\IAut(A)$ of inertial automorphisms of an abelian group $A$. They follow from Theorem 1 of \cite{DR1} and Proposition 2.2 and Theorem A of \cite{DR2}. Here we deal with finitely generated subgroups of $\IAut(A)$, while for the structure of the whole group $\IAut(A)$ see \cite{DR3}. 

\medskip

To handle the inertial automorphisms mentioned in the introduction, we 
 introduce some terminology. If $\g$ is the automorphism of the torsion-free abelian group $A$ such that $(a^\g)^{n}=a^{m}\ \forall a\in A$ (with $m,n$ coprime integers) we say that $\g$  is a {\em rational-power automorphism}. By abuse of notation, we 
write $\g=m/n$. Note that  in \cite{DR2} and \cite{DR4} we call such a $\g$ a \lq\lq multiplication" since we are regarding $A$ as a $\Q^\pi$-module, where $\Q^\pi$ is the {\it ring of rational numbers whose denominator is a $\pi$-number} and  $\pi$ is the set of primes $p$ such that $A^p=A$. 
 As in \cite{LR}, we will denote by $\Q_\pi$ the  {\it additive  group} of the ring $\Q^\pi$.

\begin{Prop}\label{Recall} {Let  $\G$ be a finitely generated group of inertial automorphisms of an abelian group $A$. Then, \\  
\noindent 
$1)$ if $A$ is torsion-free and $\pm1\ne\g\in \Aut(A)$, then $\g$ is inertial if and only if  $\g=m/n\in\Q$ and   $A^n=A^m=A$ has finite rank, thus $\IAut(A)$ is abelian in this case;\\ 
\noindent 
$2)$ if $A$ is bounded, then there is a finite $\G$-invariant subgroup $F\le A$ such that  $\G$ acts on $A/F$ by means of power automorphisms;\\ 
$3)$ if $A$ is torsion, then for each $X\le A$ there is a $\G$-invariant subgroup  $X^\G\ge X$ such that $X^\G/X$ is finite;\\ \smallskip
$4)$ if $A$ is any abelian group, then there is a $\G$-invariant torsion-free subgroup $V\le A$ such that $A/V$ is torsion.}
\end{Prop}


Clearly {\it the class of \TI-groups is closed with respect to the formation of normal subgroups and homomorphic images}. Let us give instances of \TI-groups.


\begin{Lemma}\label{esempiTtilde}
Let $G_{1}\le G_{0}$ be normal subgroups of a group $G$ with $G_{1}$ and $G/G_{0}$ finite. If any subnormal subgroup of $G_{0}/G_{1}$ is inert in $G/G_{1}$, then $G$ is \TI.
\end{Lemma}

\pf Let $H$ be a subnormal subgroup of $G$. On the one hand, $HG_1\cap G_0$ is subnormal in $G$ hence inert. On the other hand, $H$ and $HG_1\cap G_0$ are commensurable. Thus $H$ is inert, since it is commensurable with an inert subgroup.\QED



\noindent{\bf Proof of Theorem $\bf \tilde A$.} 
Suppose that $G$ is a    \TI-group.
By Theorem A of  \cite{R}, any torsion-free nilpotent normal subgroup of  $G$ is abelian. Thus $A:=\Fit (G)$ is abelian and by Proposition \ref {Recall}(1) it follows that $G/C_{G}(A)$ is abelian too.
Suppose, by the way of contradiction, that $A\ne C:=C_{G}(A)$. 
Since $G$ is hyper-(abelian or finite), there exists a $G$-invariant
subgroup $U$ of $C$ properly containing $A$ and such that $U/A$ is finite or
abelian. In the latter case $U$ is nilpotent and so
$U=A$, a contradiction. Then $U/A$ is finite, so $U/Z(U)$ is
finite and $U'$ is finite. Then $U'=1$, a contradiction again.
Hence $A=C$ and $G$ is semidihedral on $A$.

Conversely, let $G$ be a semidihedral group on a torsion-free abelian subgroup $A$. If $A$ has infinite rank, then $G$ is dihedral and every subgroup is inert. Then assume that $A$ has finite rank. Let  $H$ be a subnormal subgroup of $G$. 
 If $H\le A$, then $H$ is  inert since $G$ acts on $A$ by means of inertial automorphisms. Otherwise, by Proposition \ref{Recall}(1),  there is an element $h\in H\setminus A$ acting on $A$ as the rational power automorphism $m/n\ne 1$. If $H$ has subnormality defect $i$, we have $H\ge[H,_i A]\ge A^{(m-n)^i}$. Since $A$ has finite rank, then $A/A^{(m-n)^{i}}$ is finite and so $|HA:H|$ is finite. Then $|H^{G}:H|$ is finite and $H$ is strongly inert, hence inert. Thus $G$ is a \TI-group. 
 \QED




\section{Finitely generated groups and proof of Theorem $\bf \tilde B$}


Recall that $\tau (G)$ denotes the {\it maximum normal torsion subgroup} of a group $G$.

\begin{Prop}\label{Prop2SuSemidiedrali} Consider the following properties for a group $G$:\\
 $i)$ $G$ has a semidihedral normal subgroup with finite index $G_0$ such that $G$ acts by means of rational-power automorphisms on $A_{0}:=Fit(G_{0})$ (therefore $G$ acts trivially on $G_{0}/A_{0}$);\\
$ii)$\ $G$ has a finite normal subgroup $F$ such that $G/F$ is semidihedral. 
\medskip\\
Then $(i)$ implies $(ii)$. Moreover $(i)$ and $(ii)$ are equivalent, provided $G$ has finite rank. 
\end{Prop} 

\pf  Let $(i)$ hold. By Proposition \ref{Recall}(1), $G/C_G(A_0)$ is abelian. Thus for any $g\in G$ and $g_{0}\in G_{0}$, we have  $[g,g_{0}]\in C_{G_{0}}(A_{0})=A_{0}$. Hence $g$ acts trivially on $G_{0}/A_{0}$. 
Let $C:=C_G(A_0)$. Since $C\cap G_0=A_0$, we have that $C/A_0$ is finite. It follows that $C'$ and $F/C':=\tau(C/C')$ are
finite as well. Thus $F$ is finite. 

We claim that  {\it $\bar G:=G/F$ is semidihedral on $\bar {C}$} (we use bar notation to denote subgroups and elements modulo $F$). To show this consider two cases. If $A_0$ has finite rank, then  $G$ acts  on $A_{0}$ by means of inertial automorphisms. Otherwise, since $A_0$ has finite index in $G$, then $G$ acts on $A_0$ by means of periodic rational-power automorphisms, that is by $\pm1$. In both cases  $G$ acts on $A_{0}$ by means of inertial automorphisms. Since $C/A_{0}$ is finite,  by the same argument as in Lemma \ref{esempiTtilde}, we have that 
$G$ acts  by means of inertial automorphisms on the whole $C$. On the other hand $C_{\bar G}(\bar C)=\bar C$ as
if $\bar x\in C_{\bar G}(\bar {C})$, then $[x,A_0]\le A_0\cap F=1$. Therefore the claim follows. 

\medskip

Assume now that $(ii)$ holds and $G$ has finite rank.
Let  $n:=|F|$, $A_1/F:=Fit(G/F)$ and $C:=C_{A_{1}}(F)$. Then $F\cap C\le Z(C)$ and $C/(F\cap C)$ is abelian. Thus $[C^{n},C^{n}]\le (F\cap C)^{n^{2}}=1$. Therefore $C^{n}$ is abelian and  $A_{0}:=C^{n^{2}}$ is torsion-free abelian and has finite index, say $s$, in $A_{1}$.
By using bar notation in  $\bar G=G/A_{0}$, let $\bar G_{1}:=C_{\bar G}(\bar A_{1})$. Then $[\bar G_{1}^{s},\bar G_{1}^{s}]\le \bar A_{1}^{s^{2}}=1$ and so $\bar G_{1}^{s}$ is abelian. Moreover  $\bar G_{0}:=\bar G_{1}^{2s^{2}}$ is torsion-free abelian and has finite index in $\bar G$, since $\bar G$ has finite rank. 

Since $G_0\cap F= A_{0}\cap F=1$, then every element of $G_{0}\setminus A_{0}$ acts  fixed-point-free on $A_0\iso_G A_0F/F$. Hence $C_{G_{0}}(A_{0})=A_{0}$ and  $G_{0}$ is semidihedral on $A_{0}$.
\QED

Note that if $G$ is a group such that $G'$ has prime order, $G/G'$ is free abelian with infinite rank and $G/Z(G)$ is infinite, then $G$ has $(ii)$ but not $(i)$.


\bigskip

We now consider finitely generated groups. 

\begin{Lemma}\label{LemmaFinGen} Let $G=\la g_1,...,g_r\ra$ be a finitely generated group with 
 a torsion-free abelian normal subgroup $A$ such that $G/A$ is abelian. If  for each $i$ the element  $g_i$ acts on $A$ by means of the rational-power automorphism $m_i/n_i\in\Q$,    then $A$ is a free $\Q^\pi$-module of finite rank where $\pi=\pi(m_1...m_rn_1...n_r)$.
\end{Lemma} 

\pf  In the natural embedding of $\bar G:=G/C_G(A)$ in $\IAut(A)$, each generator $\bar g_i$
 corresponds to the rational-power automorphisms ${m_i}/{n_i}\in\Q$ (see Proposition \ref{Recall}(1)).
Thus the subring of $\rm End(A)$ generated by the image of $\bar G$ is isomorphic to $\Q^\pi$. 

 Since $G/A$ is finitely presented, we have that $A$ is finitely $G$-generated and 
$A$ is a finitely generated  $\Q^{\pi}$-module. Then $A$ is isomorphic to a direct sum of finitely many quotient of  the additive group $\Q_{\pi}$. Moreover $A$ is a free $\Q^{\pi}$-module, since it is torsion-free as an abelian group.\QED

Notice that a finitely generated semidihedral group may be obtained by a sequence of {\it finitely many HNN-extensions} starting with a  free abelian group of finite rank as a base group. However, generally such extensions are not {\it ascending} and $A$ is not finitely presented, an easy example being the extension of $\Q_{\{2,3\}}$ by the (inertial) rational-power automorphism $\g=2/3$, see Proposition 11.4.3 of \cite{LR}. On the other hand, since {finitely generated semidihedral groups} have finite rank, for such groups conditions $(i)$ and $(ii)$ of Proposition \ref{Prop2SuSemidiedrali} are equivalent. 

\vskip-0.2cm
\begin{Lemma}\label{N'finito} Let
$G$ be a finitely generated group with a normal subgroup $N$ such that $G/N$ is abelian and $G$ acts on $N/N'$  by means of inertial automorphisms. If $N'$ is finite, then  $\tau(G)$ is finite.
\end{Lemma}

\pf By arguing mod $N'$ we may assume that $N$ is abelian.

 If $N$ is torsion, then it is bounded, since it is $G$-finitely generated. Thus by Proposition \ref{Recall}(2) there is a finite $G$-subgroup $F\le N$ such that $G$ acts by means of power automorphisms on $N/F$. We may assume $F=1$. Then $G/C_G(N)$ is finite, as a group of power automorphisms of a bounded abelian group.  Therefore  the nilpotent group $C_G(N)$ is finitely generated. Thus $G$ is polycyclic and $\tau(G)$ is finite.

If $N$ is any abelian group, then, by Proposition \ref{Recall}(4), there is a torsion-free $G$-subgroup $V\le N$ such that $N/V$ is torsion. By the above $\tau (G/V)$ is finite, whence $\tau(G)$ is finite.
\QED

\vskip-0.4cm

 \begin{Lemma}\label{N2perabeliano} Let $G$ be a finitely generated \TI-group. If $G'$ is nilpotent of class 2 and $G''$ is a $p$-group, then  $\tau(G)$ is finite.
\end{Lemma}

 \pf
By Lemma \ref{N'finito}, we have that $\tau (G/G'')$ is finite. Then by Theorem $\tilde A$ and Proposition \ref{Prop2SuSemidiedrali}, there is a subgroup $G_{0}$ of finite index in $G$ such that $G_0\ge G''$ and $G_{0}/G''$ is torsion-free semidihedral. Then $T:= \tau(G_{0})\le G'' $  is  abelian. Applying again   Lemma \ref{N'finito} to $G_0/G_0''$, we have that $T/G_{0}''$ is finite. Since it suffices to prove that 
$\tau(G_0)$ is finite, we may replace $G$ by $G_{0}$. Therefore we reduce to the case where {$G/T$ is torsion-free semidihedral, where $T=\tau(G)$ is abelian and $T/G''$ is finite}. 

Then there is a finite subgroup $F$, such that $T=FG''$.  By Proposition \ref{Recall}(3) we may assume $F$ to be $G$-invariant and factor out by $F$. Thus, denoting $N:=G'$,  we have the following:\\
-  $G=\la g_1,...,g_r\ra$ is finitely generated with a nilpotent subgroup $N$ with class $2$ such that $G/N$ is abelian, $N/N'$ is torsion-free  and $N'=\tau(G)$ is a $p$-group;\\ 
- each $g_i$ acts on $N/N'$ by means of a rational-power automorphism, say $m_i/n_i\in\Q$ (by Proposition \ref{Recall}(1)).

By Lemma \ref{LemmaFinGen}, {$N/N'$ is a free $\Q^{\pi}$-module of finite rank } where  \emph{$\pi:=\pi(m_1...m_rn_1...n_r)$}.  Further, we have that  for each $a,b\in N$  and $g\in G$, if $a^{ng}=a^{m}z_{1}$ and  $b^{ng}=b^{m}z_{2}$ with $z_{1}, z_{2}\in N'\le Z(N)$ and $m,n\in\Z$, then  $[a,b]^{n^{2}g}=[a^{ng},b^{ng}]=[a^{m}z_{1}, b^{m}z_{2}]=[a^{m}, b^{m}]=[a,b]^{m^{2}}$.  Since $N'=\{[a,b]\,|,a,b\in N\}$, we have $p\not\in  \pi$.

By a classical argument (see 5.2.5 in \cite{R}), we have that {the $p$-group $\tau(G)=N'$ is finite} since it is  isomorphic to an epimorphic image of $N/N'\otimes N/N'$ which is a direct sum of finitely many copies of $\Q_{\pi}$ where $p\not\in \pi$. 
\QED

\medskip

\noindent  \textbf{Proof of Theorem $\bf \tilde B$.} 
Clearly, $(ii)$ implies $(i)$ by Lemma \ref{esempiTtilde} and Theorem $\tilde A$. 

Let $G$ be a \TI-group and note that $(ii)$ is equivalent to saying that $\tau(G)$ is finite. Clearly it suffices to prove the statement in the case where $G$ is soluble.
By induction on the derived length of $G$ we may assume that $G$ has a normal abelian subgroup $A$ such that $\tau(G/A)$ is finite. By Theorem $ \tilde A$,  Proposition \ref{Prop2SuSemidiedrali} and Lemma \ref{LemmaFinGen}, there is a subgroup $G_{0}$ of finite index in $G$ such that $G_{0}/A$ is {torsion-free semidihedral of  finite rank}. We may assume $G:=G_0$. 

Consider first the case where $A$ is a $p$-group. If $A$ is unbounded and $B\le A$ is such that $A/B$ is a Pr\"ufer group, then, by Proposition \ref{Recall}(3), such is $A/B^G$. Using bar notation in  $\bar G:=G/B^G$, we have that  $\bar G/C_{\bar G}(\bar A)$ is abelian and $\bar G'$ is nilpotent of class $2$, since  $G''\le A$. 
Then we may apply Lemma \ref{N2perabeliano} and we have that $\tau(\bar G)$ is finite. Thus $\bar A=\tau(\bar G)$ is finite, a contradiction.  Thus $A$ is bounded and, by Proposition \ref{Recall}(2), there is  a finite $G$-invariant subgroup $F\le A$ such that $G$ acts on $A/F$ by means of power automorphisms. Then $ G/C_{G}(A/F)$ is abelian and $A=\tau(G)$ is finite, again by  Lemma \ref{N2perabeliano}.

Let $A$ be {any torsion abelian group}. By the above, its primary components are finite. 
By the way of  contradiction, if $A$ is infinite, then  there is $B\le A$ such that $A/B$ is infinite with cyclic primary components. By Proposition \ref{Recall}(3), we may assume $B:=B^G$ to be $G$-invariant. Then $A/B$ has finite Pr\"ufer rank. Hence  $\bar G:=G/B$ has finite Pr\"ufer rank (and is finitely generated). Thus, by Corollary 10.5.3 of \cite{LR}, $\bar G$ is a minimax group. Therefore its torsion  normal subgroups have min (the minimal condition on subgroups) while $\bar A$ has not,  a contradiction. 

In the {general case}, since $G$ acts by means of inertial automorphisms on  $A$, by Proposition \ref{Recall}(4)  there is a torsion-free $G$-invariant subgroup $V$ of $A$ such that $A/V$ is torsion.  By the above,  $\tau (G/V)$ is finite, whence  $\tau (G)$ is finite. \QED


\end{document}